\documentclass{amsart}
\usepackage{amssymb}
\usepackage{amsfonts}
\usepackage{amsmath}
\usepackage[legalpaper,bookmarks=true,colorlinks=true,linkcolor=blue,citecolor=blue]{hyperref}
\usepackage{graphicx}%
\usepackage{fancyhdr}
\usepackage{color}
\usepackage[mathlines]{lineno}
\usepackage{lscape}
\usepackage{epsfig}
\usepackage{natbib}
\usepackage{geometry}
\usepackage{tgbonum}
\fontfamily{qcr}\selectfont

\theoremstyle{plain}

\newtheorem{proposition}{Proposition}

\numberwithin{equation}{section}

\newcommand{\Bin}{\bigskip \noindent}

\newcommand{\Ni}{\noindent}

\begin{document}
\Large
\title[Riemann and Lebesgue integrals with discrete distribution function]{Comparison between Riemann-Stieltjes and Lebesgue-Stieltjes integration using discrete distribution functions}

\begin{abstract} Integrating with respect to functions which are constant on intervals whose bounds are discontinuity points (of those functions) is frequent in many branches of Mathematics, specially in stochastic processes. For such functions and alike extension, a comparison between Riemann-Stieltjes and Lebesgue-Stieltjes integration and the integrals formulas leads to interesting facts for students (as complements of Measure Theory and Integrations) and for practitioners and and researchers. We undergone conditions of existence the Riemann-Stieltjes integrals related to that type of function and compare the results with what should be expected for Lebesgue-Stieltjes theory.\\  

\noindent $^{\dag}$ Gane Samb Lo.\\
LERSTAD, Gaston Berger University, Saint-Louis, S\'en\'egal (main affiliation).\newline
LSTA, Pierre and Marie Curie University, Paris VI, France.\newline
AUST - African University of Sciences and Technology, Abuja, Nigeria\\
gane-samb.lo@edu.ugb.sn, gslo@aust.edu.ng, ganesamblo@ganesamblo.net\\
Permanent address : 1178 Evanston Dr NW T3P 0J9,Calgary, Alberta, Canada.\\

\noindent Aladji Babacar Niang\\
LERSTAD, Gaston Berger University, Saint-Louis, S\'en\'egal.\newline
Email: aladjibacar93@gmail.com \\

\noindent Cherif Mamadou Moctar Traor\'e\\
LERSTAD, Gaston Berger University (UGB), Saint-Louis, Senegal.\\
Email : traore.cherif-mamadou-moctar@ugb.edu.sn, cheriftraore75@yahoo.com.\\

\end{abstract}
\maketitle

\section{Introduction} 

\Ni Let $(x_h)_{(h\in \mathbb{Z})}$ be an ordered sequence in $\mathbb{R}$, represented as :

$$
\cdots < x_{-(h+1)} < x_{-h} < \cdots < x_{-2} < x_{-1}<x_{0}<x_{1} < x_{2} < \cdots < x_{h} < x_{h+1}< \cdots
$$

\Bin where $h>0$. Let $F: \mathbb{R}\rightarrow \mathbb{R}$ be a continuous  over the intervals $]x_{h-1}, \ x_h[$, $h\in \mathbb{Z}$, and discontinuous at the  $x_h$'s. A special case concerns functions $F$ constant on intervals $]x_{h-1}, \ x_h[$, $h \in \mathbb{Z}$.\\

\Ni In this note, we are going to study Riemann-Stieltjes and Lebesgue-Stieljes integration with respect to such function $F$, set integration formulas and finally make comparison. This text can be seen as a complement to \cite{ips-mestuto-ang} where a general comparison is made. But the computation questions are not addressed there and we can need them is a number of disciplines, one of them being the stochastic processes with discontinuity jumps.\\

\Ni As usual, the Lebesgue-Stieltjes integration theory is simpler and more beautiful, but the Riemann-Stieltjes approach is more practical for computation purposes.\\

\Ni Then we begin with the Lebesgue-Stieltjes integration which is based on identification method and use results from \cite{ips-mestuto-ang} and books alike.

\section{Lebesgue-Stieltjes integration with a discrete distribution function}

\Ni The non-negative $F$-Lebesgue-Stieltjes measure requires that $F$ assigns to intervals $]a,b]$, $a<b$, the non-negative length $\Delta F(a,b)=F(b)-F(a)\geq 0$, meaning that $F$ is increasing in the specific case of dimension one and that $F$ is right-continuous when we use the semi-algebra $\{]a,b], \ a<b\}$ for construction the Lebesgue-Stieltjes integral (see Chapter 11, page 683, in \cite{ips-mestuto-ang}). So, here, we consider a sequence $(p_h)_{h\in \mathbb{Z}}$ of positive numbers such that\\

\Ni \textbf{A[i]} $\forall h \in \mathbb{Z}$, 

$$
p_h^{\star}=\sum_{j\leq h} p_j.<+\infty,
$$

\Ni \textbf{A[ii]} 

$$
p^{\star}=\sum_{h\in \mathbb{Z}} p_h \in ]0,+\infty].
$$

\Bin (The infinite value is not excluded). Next we set

\begin{equation}
F(x)=\sum_{h\in \mathbb{Z}} p_h^{\star} 1_{]x_{h-1},\ x_h]}(x), \ x \in \mathbb{R}. \label{dfCAD}
\end{equation}

\Bin The function $F$ in Formula \eqref{dfCAD} is right-continuous and hence is a distribution function in the sense precised above. For Riemann-Stieljes integration, the right-continuity is not required. We may use function defined as

\begin{eqnarray}
F(x)&=&\sum_{h\in \mathbb{Z}} p_h^{\star} 1_{]x_{h-1},\ x_h[}(x) \notag\\
&+&\sum_{h\in \mathbb{Z}} q_h^{\star} 1_{\{x_h\}}(x) , \ x \in \mathbb{R}, \label{dfGEN}
\end{eqnarray}

\Bin where $q_h^{\star} \in [F(x_h-0), \ F(x_h+0)]$ with the left-hand and right-hand limits $F(x_h\pm0)$ are computed with

$$
\sum_{h\in \mathbb{Z}} p_h^{\star} 1_{]x_{h-1},\ x_h[}(x), \ x \in \mathbb{R}\setminus \{x_h, \ h\in \mathbb{Z}\}.
$$

\Bin Of course, the sequence $(x_h)$ can be bounded and/or can be finite. In these cases, the notation are easily adapted.\\

\Bin In the rest of that section, we focus on two cases.\\

\Ni \textbf{Case A : F is a \textit{df}'s with jumps $(x_h)$ and constant between the jump points}. In that case, $F$ is exactly of the form \eqref{dfCAD}. Here, since the $p_h$'s are the jumps and $F$ is non-decreasing, we have that

$$
p_h^{\star}=\lambda_F(]-\infty,x_h])<+\infty,
$$

\noindent by definition a distribution function, so that \textbf{A[i]} holds. If \textbf{A[ii]} holds the Lebesgue-Stieltjes measure $\lambda_F$ is bounded. But we can easily see that $F$ in \eqref{dfCAD} is the cumulative distribution function \textit{cdf} of the countable and non-negative linear combination of the Dirac measures $\delta_{x_h}$ with coefficients $p_h$:

$$
\nu = \sum_{h\in \mathbb{Z}} p_h \delta_{x_h}.
$$

\Bin Since both functions take the null value at $-\infty$, we exactly have that

\begin{equation}
\lambda_F = \sum_{h\in \mathbb{Z}} p_h \delta_{x_h}. \label{SimpleLS}
\end{equation}
 
\Bin By using integration with respect to countable and non-negative linear combination of the Dirac measures (see Doc 05-03, page 357, in \cite{ips-mestuto-ang}, we get the following simple rule. Below, $f^{+}=max(f,0)$ and $f^{-}=max(-f,0)$ denote the positive and negative parts of a function $f$.

\begin{proposition} Let non-negative $F: \mathbb{R}\rightarrow \mathbb{R}$ satisfies: (i) $F$ is  right-continuous, (ii) $F$ assigns to intervals $]a,b]$, $a<b$, non-negative length $\Delta F(a,b)=F(b)-F(a)\geq 0$, (iii) $F$ is constant on intervals $]x_{h-1}, \ x_h[$ and, (iv) is discontinuous at the  $x_h$'s. Hence:\\

\Ni (a) for any constant-sign function $f :\mathbb{R}\rightarrow \overline{\mathbb{R}}$, the Lebesgue-Stieltjes integral of $f$ exists in 
$\overline{\mathbb{R}}_+$ or in $\overline{\mathbb{R}}_{-}$ and is given by

\begin{equation}
\int f(x) dF(x)= \sum_{h \in \mathbb{Z}} p_h f(x_h); \label{LSInteg}
\end{equation}

\Bin (b) for any function $f :\mathbb{R}\rightarrow \overline{\mathbb{R}}$, such that

$$
0\leq \sum_{h \in \mathbb{Z}} p_h f^{+}(x_h)<\infty \ \ or \ \ 0\leq \sum_{h \in \mathbb{Z}} p_h f^{-}(x_h)<\infty,
$$

\noindent  the Lebesgue-Stieljtes integral of $f$ exists in $\overline{\mathbb{R}}$ and is given by Formula \eqref{LSInteg}.
\end{proposition}

\Bin The function $F$ is continuous on $]x_{h-1}. x_{h}[$. Let us denote by $\lambda_{F^{(h)}}$ the Lebesgue-Stieltjes measure on $]x_{h-1}. x_{h}[$. Then Formula \eqref{SimpleLS} becomes

\begin{equation}
\lambda_F=\sum_{h\in \mathbb{Z}} \biggr(p_h \delta_{x_h}+ \lambda_{F^{(h)}}\biggr). \label{SimpleLS1}
\end{equation}

\Bin This leads to the following result.

\begin{proposition} Let $F: \mathbb{R}\rightarrow \mathbb{R}$ be non-negative such that: (i) $F$ right-continuous and let $F$ assigns to intervals $]a,b]$, $a<b$, 
non-negative length $\Delta F(a,b)=F(b)-F(a)\geq 0$, (ii) $F$ is continuous on intervals $]x_{h-1}, \ x_h[$ and (iii) $F$ is discontinuous at the  $x_h$'s. Hence:\\

\Ni (a) for any measurable function $f :\mathbb{R}\rightarrow \overline{\mathbb{R}}$, of constant sign, almost-surely continuous on each $]x_{h-1}, \ x_h[$ with respect to the Lebesgue measure, which ensures that

$$
\int 1_{]x_{h-1}, \ x_h[} f(x) \ dF(x) \in \mathbb{R},
$$

\Bin the Lebesgue-Stieljtes integral of $f$ exists in $\overline{\mathbb{R}}$ and is given by

\begin{equation}
\int f(x) dF(x)= \sum_{h \in \mathbb{Z}} \biggr( p_h f(x_h) + \int 1_{]x_{h-1}, \ x_h[} f(x) \ dF(x) \biggr); \label{LSGInteg}
\end{equation}

\Bin (b) for any function measurable $f :\mathbb{R}\rightarrow \overline{\mathbb{R}}$, almost surely continuous on each $]x_{h-1}, \ x_h[$ with respect to the Lebesgue measure, which ensures that

$$
\int 1_{]x_{h-1}, \ x_h[} f(x) \ dF(x) \in \mathbb{R},
$$

\Bin such that

$$
0\leq \sum_{h \in \mathbb{Z}} \biggr( p_h f^{+}(x_h) + \int 1_{]x_{h-1}, \ x_h[} f^{+}(x) \ dF(x) \biggr) \leq +\infty
$$

\Bin or

$$
\sum_{h \in \mathbb{Z}} \biggr( p_h f^{-}(x_h) + \int 1_{]x_{h-1}, \ x_h[} f^{-} \ dF(x) \biggr),
$$

\noindent  then the Lebesgue-Stieltjes integral of $f$ exists in $\overline{\mathbb{R}}$ and is given by Formula \eqref{LSGInteg}.
\end{proposition}

\Bin Let us now move to the Riemann-Stieltjes integral.

\section{Riemann-Stieltjes integral}

\Ni In the frame of Riemann or Riemann-Stieltjes integration, the domain of integration is bounded and is taken as $[a,b]$, $a<b$ and the function to be integrated is bounded, say by $M$. So we have

$$
f: [a,b] \rightarrow \mathbb{R}, \ \ |f|\leq M.
$$

\Bin We consider $F$ not necessarily increasing, having a finite number $p$ of discontinuity $x_1<x_2<\cdots<x_p$, $p\geq 1$. We denote $x_0=a$ and $b=x_{p+1}$. Although the Riemann-Stieltjes is more restrictive, we will see that we do not need $F$ to be right or left continuous on the $x_h$. Further, we know that $f$ is integrable if only if $f$ is continuous almost surely when $F$ is continuous. Here we will see that, if $F$ is constant over the intervals $]x_{h-1}, x_{h}[$, only the continuity of $f$ on the $x_h$'s is taken as sufficient.\\

\Ni Let us consider arbitrary subdivisions

$$
\pi_n=a=(y_{0,n}<\cdots<y_{j,n}<\cdots < y_{\ell(n)-1,n}<y_{\ell(n),n}=b), \ n\geq 1
$$ 

\Bin and choose $c_n=(c_{j,n})_{0\leq j \leq \ell(n)-1}$ with $c_{j,n} \in ]y_{j,n},y_{j+1,n}[$. We require that $\ell(n)$ and the sequence of modulii of the subdivisions both converge to zero, that is, as $n\rightarrow +\infty$

$$
\ell(n)\rightarrow +\infty \ and \ p(\pi_n)=\max_{0\leq j \leq \ell(n)-1} (y_{j+1,n}-y_{j,n}) \rightarrow 0. 
$$

\Ni Let us study the Riemann-Stieltjes sums

\begin{equation}
S(\pi_n,c_n,f,F)=\sum_{0\leq j \leq \ell(n)-1} f(c_{j,n}) (F(y_{j+1,n})-F(y_{j,n})). \label{Rsum}
\end{equation}

\Bin We begin with:\\

\Ni Case 1. $F$ is constant over the intervals $]x_{h-1}, x_{h}[$. To make things simpler, we provisionally require that $F(a)=F(a+0)$ and $F(b)=F(b-0)$. The idea we exploit is that, if an interval $]y_{j,n},y_{j+1,n}[$ is in one of the $[x_{h-1}, \ x_h[$, the corresponding term which is

\begin{equation}
f(c_{j,n}) (F(y_{j+1,n})-F(y_{j,n})). \label{RsumElem}
\end{equation}

\Bin vanishes. So the Riemann-Stieltjes reduces to the sum of terms \eqref{RsumElem} for which :\\

\noindent (i) either $y_{j,n}=x_h$  for some $h \in \{1,\cdots,p\}$. [There are no problems at $x_0$ or at $x_{p+1}$ with the assumption at those points],\\

\noindent (ii) or $x_h \in ]y_{j,n},y_{j+1,n}[$ for some $h \in \{x_1,\cdots,x_p\}$.\\

\noindent We are going to see what can happen around each $h \in \{1,\cdots,p\}$. In the case (i), the two remaining terms are

$$
L_1(n)=\{F(x_h)-F(y_{j,n})\} f(c_{j,n}) + \{F(y_{j+2,n})-F(x_h)\} f(c_{{j+1},n}).
$$ 

\Bin we study one of those two terms and translate the result to the other. Let us call the first as $L_1(1,n)$. By taking all $c_{j,n}$ in that case as $c_{j,n}=x_h$, we have that $L_1(1,n)$ converges to

$$
\{F(x_h)-F(x_h-0)\} f(x_h).
$$

\Bin By choosing $c_{j,n}<x_h$ for each subdivision $\pi_n$, and by supposing that $f$ is left-continuous at $x_h$, wet get that $L_1(1,n)$ converges to 

$$
\{F(x_h)-F(x_h-0)\} f(x_h-0).
$$

\Bin In summary for the different choices, by doing the same for the second term of $L_1(n)$ , and by assuming that $f$ has left-limits and right-limits at the $x_h$, $L_1(1,n)$ converges to the four possible values

\begin{eqnarray*}
&&A + B,\\
&&  A \in \biggr\{\{F(x_h)-F(x_h-0)\} f(x_h),  \ \{F(x_h)-F(x_h-0)\} f(x_h-0)\biggr\},\\
&& B \in \biggr\{\{F(x_h+0)-F(x_h)\} f(x_h),  \ \{F(x_h+0)-F(x_h)\} f(x_h+0)\biggr\}
\end{eqnarray*}

\Bin In the case (ii), we have three choices of $c_{j,n}$ we have for all subdivisions.\\

\Ni For $c_{j,n}=x_h$, the remaining term \eqref{RsumElem} goes to

$$
\{F(x_h+0)-F(x_h-0)\} f(x_h)
$$

\Bin For $c_{j,n} \in [y_{j-1,n},\ x_h[$, the remaining term \eqref{RsumElem} goes to

$$
\{F(x_h+0)-F(x_h-0)\} f(x_h-0)
$$

 For $c_{j,n} \in ]x_h,y_{j,n}]$
$$
\{F(x_h+0)-F(x_h-0)\} f(x_h+0)
$$

\Bin In summary, if $f$ is continuous at the $x_h$, all the remaining terms converge to

$$
\{F(x_h+0)-F(x_h-0)\} f(x_h)
$$ 

\Bin We may write our conclusion where we lift the regularity conditions of $F$ on $a$ and $b$.

\begin{proposition} Let $a<b$, $f$ be bounded on $[a,b]$, $x_1<\cdots<x_p$ be $p$ points in $]a,b[$. Let $F$ be a real-valued function defined on $[a,b]$, discontinuous only on the $x_h$'s with discontinuity jumps $p_h=\{F(x_h+0)-F(x_h-0)\}$ and constant on the intervals $]x_{h-1},x_{h}[$ with $x_0=a$ and $x_{p+1}=b$, If $f$ is continuous on any of the $x_h$, $h \in \{1,\cdots,p\}$, then any sequence of Riemann-Stieltjes sums with unbounded number of elements ($\ell(n)\rightarrow +\infty$) and modulii going to zero, converges to

\begin{eqnarray}
I(f,F)&=& \{F(a+0)-F(a)\}f(a) + \{F(b)-F(b-0)\}f(b)\\
&+& \sum_{h=1}^{p} p_h f(x_h) \notag
\end{eqnarray}

\noindent which, by definition, is the Riemann-Stieltjes integral :

$$
I(f,F)=\int_{a}^{b} f(x) \ dF(x).
$$

\noindent If we change the constancy of $F$ on the $]x_{h-1},x_{h}[$ to continuity simply, we have that whenever $f$ is continuous at the $x_h$'s and almost-surely continuous on each $]x_{h-1}, \ x_h[$, we get

\begin{eqnarray}
\int_{a}^{b} f(x) \ dF(x)&=& \{F(a+0)-F(a)\}f(a) + \{F(b)-F(b-0)\}f(b)\\
&+& \sum_{h=1}^{p} p_h f(x_h) \notag \\
&+& \sum_{h=1}^{p+1} \int_{x_{h-1}}^{x_h} f(x) \ dF(x). \notag
\end{eqnarray}
\end{proposition} 

\Bin We may move to an improper Riemann-Stieltjes integral by extending the discontinuity points $(x_h)$ to an infinite on $(x_h)_{h\in \mathbb{Z}}$ such that

$$
\forall a<b, \# \biggr(\{x_h, \ \ h \in \mathbb{Z}\} \bigcap [a,b]\biggr) <+\infty.
$$ 

\noindent If $f$ is continuous at the $x_h$ and almost-surely continuous on each $]x_{h-1}, \ x_h[$, all the integrals

$$
\int_{a}^{b} f(x) \ dF(x)
$$

\Bin do exist. By choosing $a$ and $b$ outside the points $x_h$, the improper integral, taken as the limit when $a\downarrow -\infty$ and $b\uparrow +\infty$, is

$$
\int_{-\infty}^{+\infty} f(x) \ dF(x),
$$

\Bin and we have 
\begin{equation}
\int f(x) dF(x)= \sum_{h \in \mathbb{Z}} \biggr( p_h f(x_h) + \int_{x_{h-1}}^{x_h} f(x) \ dF(x) \biggr). \label{LSGRSInteg}
\end{equation}

\section{Some conclusion} 

\Ni We draw a number of facts from the previous developments.\\

\Ni (a) As said earlier, the Lebesgue-Stieltjes is more comfortable at the price of $F$ is \textit{distribution function}, right-continuous. How, ever if $F$ is not, the Riemann-Stieltjes offers more chances. But in turn, the Riemann-Stieltjes is not easy to handle for an infinite number of jumps on compact sets whereas Lebesgue-Stieltjes does not case as long as integrability conditions holds.\\

\Ni (b) It is remarkable that for a function $F$ having finite number of jump points in a compact set $[a,b]$, $a<b$, and constant between jump points, the existence of the Riemann-Stieltjes integral requires only the continuity on the jump points, regardless of the monotonicity of $F$ and of the left or right continuity of $F$ at the jump points.\\

\Ni (c) The Lebesgues-Stieltjes integral \eqref{LSGInteg} on $\mathbb{R}$ and the improper Riemann-Stieltjes \eqref{LSGRSInteg} coincide under the integrability conditions of the Lebesgues-Stieltjes.\\

\end{document}